\documentclass[a4paper]{article}

\usepackage{xcolor}
\definecolor{kth-blue}{RGB/cmyk}{25,84,166/0.849,0.494,0,0.349}
\definecolor{kth-red}{RGB/cmyk}{157,16,45/0,0.898,0.713,0.384}
\definecolor{kth-green}{RGB/cmyk}{98,146,46/0.329,0,0.685,0.427}
\definecolor{kth-lightblue}{RGB/cmyk}{36,160,216/0.833,0.259,0,0.153}
\definecolor{kth-lightred}{RGB/cmyk}{228,54,62/0,0.763,0.728,0.106}
\definecolor{kth-lightgreen}{RGB/cmyk}{176,201,43/0.124,0,0.786,0.212}
\definecolor{kth-pink}{RGB/cmyk}{216,84,151/10,0.611,0.301,0.153}
\definecolor{kth-yellow}{RGB/cmyk}{250,185,25/0,0.26,0.9,0.0196}
\definecolor{kth-darkgray}{RGB/cmyk}{101,101,108/0.0648,0.0648,0,0.576}
\definecolor{kth-middlegray}{RGB/cmyk}{189,188,188/0,0.00529,0.00529,0.259}
\definecolor{kth-lightgray}{RGB/cmyk}{227,229,227/0.00873,0,0.00873,0.102}
\colorlet{blue}{kth-blue}
\colorlet{red}{kth-red}
\colorlet{green}{kth-green}
\colorlet{lightblue}{kth-lightblue}
\colorlet{lightred}{kth-lightred}
\colorlet{lightgreen}{kth-lightgreen}
\colorlet{pink}{kth-pink}
\colorlet{yellow}{kth-yellow}
\colorlet{gray}{kth-darkgray}
\colorlet{middlegray}{kth-middlegray}
\colorlet{lightgray}{kth-lightgray}

\usepackage{custompkg}
\usepackage{customdef}
\usepackage{authblk}
\usepackage{changes}
\definechangesauthor[color = blue, name = {Mikael}]{MJ}
\definechangesauthor[color = green, name = {Arda}]{AA}
\definechangesauthor[color = red, name = {Martin}]{MB}

\SetAlgoLined
\RestyleAlgo{ruled}

\title{{POLO: a {POL}icy-based {O}ptimization library}}
\author[1]{Arda Aytekin\thanks{Email:
  \href{mailto:aytekin@kth.se}{aytekin@kth.se}, corresponding author.}}
\author[1]{Martin Biel\thanks{Email:
  \href{mailto:mbiel@kth.se}{mbiel@kth.se}.}}
\author[1]{Mikael Johansson\thanks{Email:
  \href{mailto:mikaelj@kth.se}{mikaelj@kth.se}.}}
\affil[1]{
KTH Royal Institute of Technology\\
School of Electrical Engineering and Computer Science\\
SE-100 44 Stockholm
}

\begin{document}
\maketitle

\begin{abstract}
  We present \texttt{POLO} --- a C++ library for large-scale parallel
  optimization research that emphasizes ease-of-use, flexibility and efficiency
  in algorithm design. It uses multiple inheritance and template programming to
  decompose algorithms into essential policies and facilitate code reuse. With
  its clear separation between algorithm and execution policies, it provides
  researchers with a simple and powerful platform for prototyping ideas,
  evaluating them on different parallel computing architectures and hardware
  platforms, and generating compact and efficient production code. A C-API is
  included for customization and data loading in high-level languages.
  \texttt{POLO} enables users to move seamlessly from serial to multi-threaded
  shared-memory and multi-node distributed-memory executors. We demonstrate how
  \texttt{POLO} allows users to implement state-of-the-art asynchronous parallel
  optimization algorithms in just a few lines of code and report experiment
  results from shared and distributed-memory computing architectures.  We
  provide both \texttt{POLO} and \texttt{POLO.jl}, a wrapper around
  \texttt{POLO} written in the Julia language, at
  \href{https://github.com/pologrp}{https://github.com/pologrp} under the
  permissive MIT license.
\end{abstract}

\section{Introduction}

Wide adoption of Internet of Things (IoT) enabled devices, as well as
developments in communication and data storage technologies, have made the
collection, transmission, and storage of bulks of data more accessible than
ever. Commercial cloud-computing providers, which traditionally offer their
available computing resources at data centers to customers, have also started
supporting low-power IoT-enabled devices available at the customers' site in
their cloud ecosystem~\cite{2018-Greengrass,2018-IoTEdge,2018-IoTHub}. As a
result, we are experiencing an increase in not only the problem dimensions of
data-driven machine-learning applications but also the variety of computing
architectures on which these applications are deployed.

There is no silver bullet for solving machine-learning problems. Because the
problems evolve continuously, one needs to design new algorithms that are
tailored for the specific problem structure, and exploit all available computing
resources. However, algorithm design is a non-trivial process. It consists in
prototyping new ideas, benchmarking their performance against that of the
state-of-the-art, and finally deploying the successful ones in production. Many
ideas today are prototyped in high-level languages such as MATLAB, Python and
Julia, and it is very rare that researchers implement their algorithms in
lower-level languages such as C and C++. Even so, these prototypes are often
incomplete, have different abstractions and coding style, and are hand-crafted
for specific problem-platform combinations. Ultimately, this results in either
performance degradations of high-level implementations in production, or high
engineering costs attached to rewriting the low-level implementations from
scratch for each different problem-platform combination.

Currently, there exist numerous machine-learning libraries, each targeting a
different need. Libraries such as \texttt{PyTorch/Caffe2}~\cite{2018-PyTorch},
\texttt{Theano}~\cite{2016-Team} and \texttt{TensorFlow}~\cite{2016-Abadi}
primarily target deep-learning applications, and support different powerful
computing architectures. High-level neural-network frameworks such as
\texttt{Keras}~\cite{2015-Chollet} provide user-friendly interfaces to backends
such as \texttt{Theano} and \texttt{TensorFlow}. Even though these libraries
implement many mature algorithms for solving optimization problems resulting
from deep-learning applications, they fall short when one needs to prototype and
benchmark novel solvers~\cite{2017-Curtin,2018-Vasilache}. Algorithm designers
need to either use the provided communication primitives
explicitly~\cite{2018-PyTorch} or interfere with the computation graph
facilities~\cite{2016-Team,2016-Abadi} to write the algorithms from scratch.
Recent libraries and domain-specific languages such as
\texttt{Ray}~\cite{2017-Moritz} and \texttt{Tensor
Comprehensions}~\cite{2018-Vasilache} aim at extending the capabilities of
computation graph facilities of aforementioned libraries, and targeting users'
custom needs regarding network architectures and data shapes, respectively. More
lightweight options such as \texttt{mlpack}~\cite{2017-Curtin} and
\texttt{Jensen}~\cite{2018-Iyer}, on the other hand, aim at providing more
generic optimization frameworks in C++. Unfortunately, these options are not
centered around algorithm design, either. Although both options provide a good
selection of predefined algorithms, these algorithms are implemented for
serial~\cite{2018-Iyer} or, to a limited extent, single-machine
parallel~\cite{2017-Curtin} computations. Designers still need to rewrite these
implementations for, say, multi-machine parallel computations. In short, we
believe that there is a need for a generic optimization platform that
facilitates prototyping, benchmarking and deployment of algorithms on
different platforms without much performance penalty.

\begin{figure*}
  \centering
  \includegraphics{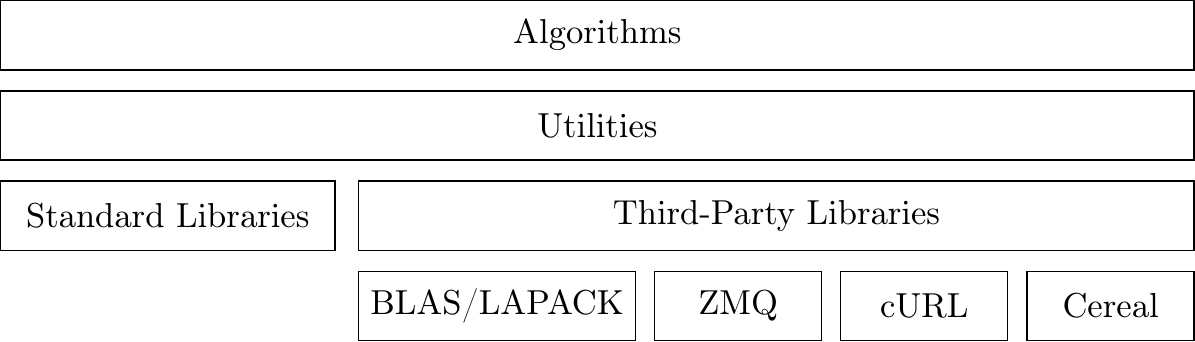}
  \caption{Structure of \texttt{POLO}. Utilities layer provides thin wrappers
  for the essential functionalities provided by the standard and third-party
  libraries, and implements custom ones, which are then used in the algorithms
  layer.}\label{fig:polo}
\end{figure*}

In this paper, we present \texttt{POLO} --- an open-source, header-only C++
library that uses the policy-based design technique~\cite{2005-Alexandrescu}. It
consists of two primary layers (Figure~\ref{fig:polo}). The utilities layer
builds on standard libraries and lightweight third-party libraries, and offers a
set of well-defined primitives for atomic floating-point operations, matrix
algebra, communication and serialization, logging and dataset reading.
Throughout the development, we have put special emphasis on making the utilities
layer portable and efficient. The algorithms layer, on the other hand, builds on
the utilities layer and implements different families of algorithms. The library
abides by modern C++ design principles~\cite{2005-Alexandrescu}, and follows
ideas from recent parallel algorithms library
proposals~\cite{2013-Hoberock,2015-Hoberock} to
\begin{itemize}
  \item decouple optimization algorithm building blocks from system
    primitives,
  \item facilitate code reuse, and,
  \item generate tight and efficient code.
\end{itemize}

In the rest of the paper, we introduce \texttt{POLO} in more detail. In
Section~\ref{sec:motivation}, we motivate the design of \texttt{POLO} based on
our observations from a family of optimization algorithms. In
Section~\ref{sec:design}, we provide detailed information about the design of
the library and the supported computing architectures. In
Section~\ref{sec:examples}, we show how some of the state-of-the-art algorithms
can be quickly prototyped in \texttt{POLO} together with their performance
comparisons against each other. Finally, in Section~\ref{sec:conclusion}, we
conclude the paper with further discussions.

\section{Motivation}\label{sec:motivation}

To demonstrate the power of policy-based design for distributed optimization
library development, we consider \emph{regularized optimization problems} on the
form:
\begin{align}
  \begin{aligned}\label{eqn:comp-opt}
    & \minimize_{x \in \set{R}^{d}}
    & & \func[\phi]{x} \coloneqq \func[F]{x} + \func[h]{x} \,.
  \end{aligned}
\end{align}
Here, $\func[F]{\cdot}$ is a differentiable function of $x$ and
$\func[h]{\cdot}$ is a possibly non-differentiable function. In machine-learning
applications, the smooth part of the loss typically represents the empirical
data loss, and is a finite sum of loss functions
\begin{align*}
  \func[F]{x} = \sum_{n = 1}^{N} \func[f_{n}]{x} \,,
\end{align*}
while the non-smooth part is a regularizer that encourages certain properties
(such as sparsity) of the optimal solution.

One approach for solving problems on the form~\eqref{eqn:comp-opt} is to use
\emph{proximal gradient methods}. The basic form of the the proximal gradient
iteration is
\begin{align}\label{eqn:prox-grad}
  x_{k+1} = \argmin[x \in \set{R}^{d}]{\func[F]{x_{k}} +
    \dotprod{\nabla\func[F]{x_{k}}}{x - x_{k}} + \func[h]{x} +
    \frac{1}{2\gamma_{k}}\norm{x-x_{k}}_{2}^{2}} \,,
\end{align}
where $\gamma_{k}$ is a step-size parameter. Thus, the next iterate, $x_{k+1}$,
is selected to be the minimizer of the sum of the first-order approximation of
the differentiable loss function around the current iterate, $x_{k}$, the
non-differentiable loss function, and a quadratic penalty on the deviation from
the current iterate~\cite{2017-Beck}. After some algebraic manipulations, one
can rewrite \eqref{eqn:prox-grad} in terms of the proximal
operator~\cite{2017-Beck}
\begin{align*}
  x_{k+1} & = \func[\operatorname{prox}_{\gamma_{k}h}]{x_{k} -
    \gamma_{k}\nabla\func[F]{x_{k}}} \\
          & \coloneqq \argmin[x \in \set{R}^{d}]{\gamma_{k}\func[h]{x_{k}} +
    \frac{1}{2}\norm{x - \left( x_{k} - \gamma_{k} \nabla\func[F]{x_{k}}
      \right)}_{2}^{2}} \,.
\end{align*}
As a result, the method can be interpreted as a two-step procedure: first, a
query point is computed by modifying the current iterate in the direction of the
negative gradient, and then the prox operator is applied to this query point.

Even though the proximal gradient method described in~\eqref{eqn:prox-grad}
looks involved, in the sense that it requires solving an optimization problem at
each iteration, the prox-mapping can actually be evaluated very efficiently for
several important functions $\func[h]{\cdot}$. Together with its strong
theoretical convergence guarantees, this makes the proximal gradient method a
favorable option in many applications. However, the gradient calculation step in
the vanilla proximal gradient method can be prohibitively expensive when the
number of data points ($N$) or the dimension of the decision vector ($d$) is
large enough. To improve the scalability of the proximal gradient method,
researchers have long tried to come up with ways of parallelizing the proximal
gradient computations and more clever query points than the simple gradient step
in~\eqref{eqn:prox-grad}~\cite{2007-Blatt,2016-Aytekin,
  2014-Defazio,1964-Polyak,1983-Nesterov,2011-Duchi,2012-Zeiler,2014-Kingma,
2016-Dozat, 2015-Sra,2011-Recht,2016-Leblond,2017-Pedregosa}. As a result, the
proximal gradient family encompasses a large variety of algorithms. We have
listed some of the more influential variants in Table~\ref{tbl:proxgrad}.

\subsection{A Look at Proximal Gradient Methods}

\begin{table}
  \centering
  \caption{Some members of the proximal gradient methods. \texttt{s},
    \texttt{cr}, \texttt{ir} and \texttt{ps} under the \texttt{execution} column
    stand for serial, consistent read/write, inconsistent read/write and
    Parameter Server~\cite{2013-Li}, respectively.}\label{tbl:proxgrad}
  \begin{tabular}{@{}llllll@{}}
    \toprule
    & \multicolumn{2}{c}{$g$} & & $\func[h]{x}$ \\
    \cmidrule(lr){2-3} \cmidrule(lr){5-5}
    Algorithm                       & \texttt{boosting}   & \texttt{smoothing} & \texttt{step}          & \texttt{prox} & \texttt{execution} \\ \midrule
    SGD                             & $\times$            & $\times$           & $\gamma$, $\gamma_{k}$ & $\times$      & \texttt{s, cr, ps}  \\
    IAG~\cite{2007-Blatt}           & \texttt{aggregated} & $\times$           & $\gamma$               & $\times$      & \texttt{s, cr, ps}  \\
    PIAG~\cite{2016-Aytekin}        & \texttt{aggregated} & $\times$           & $\gamma$               & $\checkmark$  & \texttt{s, cr, ps}  \\
    SAGA~\cite{2014-Defazio}        & \texttt{saga}       & $\times$           & $\gamma$               & $\checkmark$  & \texttt{s, cr, ps}  \\
    Momentum~\cite{1964-Polyak}     & \texttt{classical}  & $\times$           & $\gamma$               & $\times$      & \texttt{s}          \\
    Nesterov~\cite{1983-Nesterov}   & \texttt{nesterov}   & $\times$           & $\gamma$               & $\times$      & \texttt{s}          \\
    AdaGrad~\cite{2011-Duchi}       & $\times$            & \texttt{adagrad}   & $\gamma$               & $\times$      & \texttt{s}          \\
    AdaDelta~\cite{2012-Zeiler}     & $\times$            & \texttt{adadelta}  & $\gamma$               & $\times$      & \texttt{s}          \\
    Adam~\cite{2014-Kingma}         & \texttt{classical}  & \texttt{rmsprop}   & $\gamma$               & $\times$      & \texttt{s}          \\
    Nadam~\cite{2016-Dozat}         & \texttt{nesterov}   & \texttt{rmsprop}   & $\gamma$               & $\times$      & \texttt{s}          \\
    AdaDelay~\cite{2015-Sra}        & $\times$            & $\times$           & $\gamma_{k}$           & $\checkmark$  & \texttt{s, cr, ps}  \\
    HOGWILD!~\cite{2011-Recht}      & $\times$            & $\times$           & $\gamma$               & $\times$      & \texttt{ir}         \\
    ASAGA~\cite{2016-Leblond}       & \texttt{saga}       & $\times$           & $\gamma$               & $\times$      & \texttt{ir}         \\
    ProxASAGA~\cite{2017-Pedregosa} & \texttt{saga}       & $\times$           & $\gamma$               & $\checkmark$  & \texttt{ir}         \\
    \bottomrule
  \end{tabular}
\end{table}

A careful review of the serial algorithms in the proximal gradient family
reveals that they differ from each other in their choices of five distinctive
algorithm primitives: (1) which gradient surrogate they use; (2) how they
combine multiple gradient surrogates to form a search direction, a step we refer
to as \emph{boosting}; (3) how this search direction is filtered or scaled,
which we call \emph{smoothing}; (4) which step-size policy they use; and (5) the
type of projection they do in the prox step. For instance, stochastic gradient
descent (SGD) algorithms use partial gradient information coming from functions
($N$) or decision vector coordinates ($d$) as the gradient surrogate at each
iteration, whereas their aggregated
versions~\cite{2007-Blatt,2016-Aytekin,2014-Defazio} accumulate previous partial
gradient information to boost the descent direction. Similarly, different
momentum-based methods such as the classical~\cite{1964-Polyak} and
Nesterov's~\cite{1983-Nesterov} momentum accumulate the full gradient
information over iterations. Algorithms such as AdaGrad~\cite{2011-Duchi} and
AdaDelta~\cite{2012-Zeiler}, on the other hand, use the second-moment
information from the gradient and the decision vector updates to adaptively
scale, \ie, smooth, the gradient surrogates. Popular algorithms such as
Adam~\cite{2014-Kingma} and Nadam~\cite{2016-Dozat}, available in most
machine-learning libraries, incorporate both boosting and smoothing to get
better update directions. Algorithms in the serial setting can also use
different step-size policies and projections independently of the choices above,
which results in the pseudocode representation of these algorithms given in
Algorithm~\ref{alg:prox-grad}.
\begin{algorithm}[ht]
  \LinesNumbered
  \SetKwFunction{NotDone}{not\_done}
  \SetKwFunction{Gradient}{gradient\_surrogate}
  \SetKwFunction{Boost}{boosting}
  \SetKwFunction{Smooth}{smoothing}
  \SetKwFunction{Step}{step}

  \KwData{Differentiable functions, $\func[f_{n}]{\cdot}$; regularizer,
    $\func[h]{\cdot}$.}
  \KwIn{Initial decision vector, $x_{0}$; step size, $\gamma_{k}$.}
  \KwOut{Final decision vector, $x_{k}$.}

  $k \gets 0$\;
  $g \gets 0$\;
  \While{\NotDone{$k$, $g$, $x_{k}$, $\func[\phi]{x_{k}}$}}{
    $g \gets$ \Gradient{$x_{k}$}\tcp*[r]{partial or full gradient}
    $g \gets$ \Boost{$k$, $g$}\tcp*[r]{optional}\label{alg:startline}
    $g \gets$ \Smooth{$k$, $g$, $x_{k}$}\tcp*[r]{optional}
    $\gamma_{k} \gets$ \Step{$k$, $g$, $x_{k}$, $\func[\phi]{x_{k}}$}\;
    $x_{k + 1} \gets \func[\operatorname{prox}_{\gamma_{k}h}]{x_{k} -
      \gamma_{k}g}$\tcp*[r]{prox step}\label{alg:endline}
    $k \gets k + 1$\;
  }
  \Return{$x_{k}$}\;
  \caption{Serial implementation of proximal gradient methods.}
  \label{alg:prox-grad}
\end{algorithm}

Most of the serial algorithms in this family either have direct counterparts in
the parallel setting or can be extended to have the parallel \emph{execution}
support. However, adding parallel execution support brings yet another layer of
complexity to algorithm prototyping. First, there are a variety of parallel
computing architectures to consider, from shared-memory and distributed-memory
environments with multiple CPUs to distributed-memory heterogeneous computing
environments which involve both CPUs and GPUs. Second, some of these
environments, such as the shared-memory architectures, offer different choices
in how to manage race conditions. For example, some
algorithms~\cite{2007-Blatt,2016-Aytekin,2014-Defazio,2015-Sra} choose to use
mutexes to \emph{consistently} read from and write to the shared decision vector
from different processes, whereas
others~\cite{2011-Recht,2016-Leblond,2017-Pedregosa} prefer atomic operations to
allow for \emph{inconsistent} reads and writes. Finally, the choice in a
specific computing architecture might constrain choices in other algorithm
primitives. For instance, if the algorithm is to run on a distributed-memory
architecture such as the Parameter Server~\cite{2013-Li}, where only parts of
the decision vector and data are stored in individual nodes, then only updates
and projections that embrace data locality can be used.

\subsection{Algorithmic perspective}

\begin{figure}[ht]
  \centering
  \includegraphics[width = \textwidth]{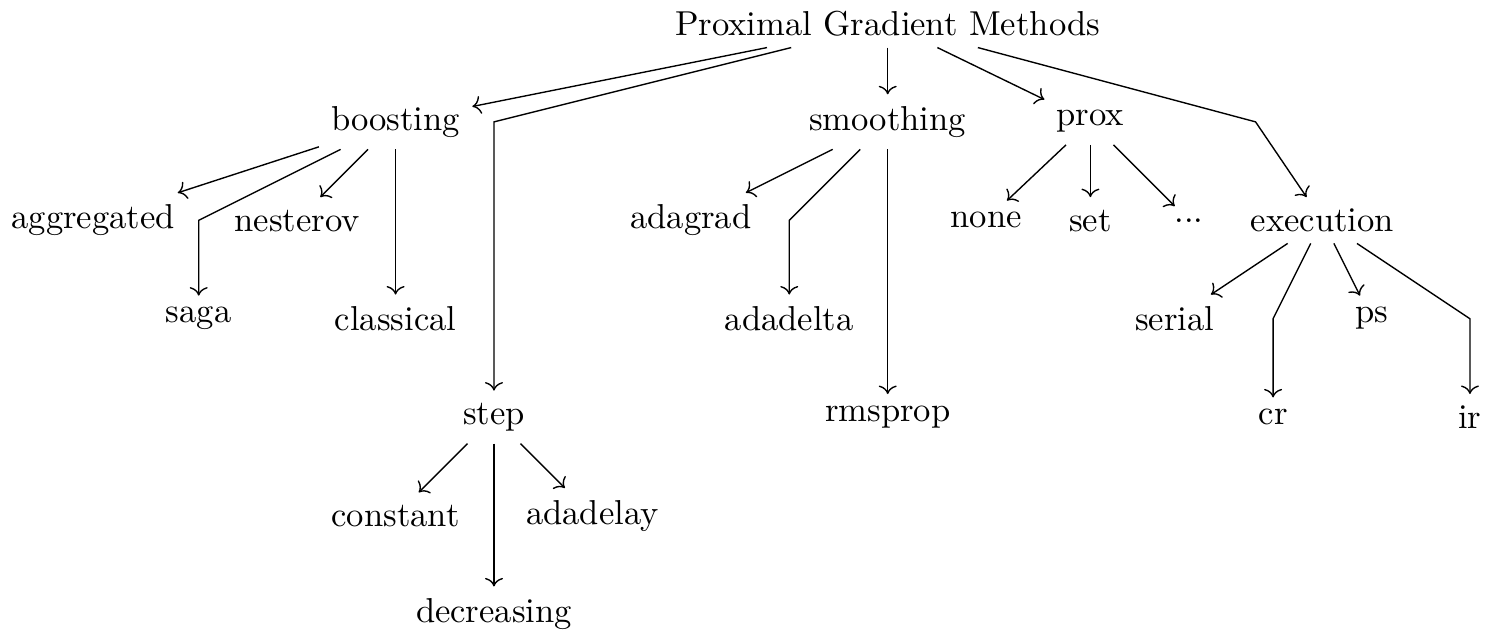}
  \caption{Tree representation of different choices in proximal gradient
    methods. Both Adam~\cite{2014-Kingma} and Nadam~\cite{2016-Dozat} use
    \texttt{rmsprop} smoothing, \texttt{constant} step and \texttt{none} prox
    policies on \texttt{serial} executors. They differ from each other in their
    respective choices of \texttt{classical} and \texttt{nesterov} boosting
    policies. Recent research~\cite{2018-Reddi} suggests that using a different
    smoothing policy, \ie, \texttt{adamax}, results in better performance
    compared to Adam-like algorithms.}\label{fig:proxgrad}
\end{figure}

Out of the large number of proximal gradient methods proposed in the literature,
only very few have open source low-level implementations, and the existing codes
are difficult to extend. This means that algorithm designers need to reimplement
others' work to benchmark their ideas even when the ideas differ from each other
only in a few aspects. Still, our observations in the previous subsection reveal
that algorithm design can be seen as a careful selection of a limited number of
compatible policies, see Figure~\ref{fig:proxgrad}.

\texttt{POLO} aims at being an optimization library centered around algorithm
design, that facilitates rapid experimentation with new ideas without
sacrificing performance. It has a clear separation between algorithm and system
primitives to allow for (1) prototyping of new ideas with as few lines of code
as possible, and (2) benchmarking resulting algorithms against each other on
different platforms as simple as possible. Ultimately, this will lead to faster
development of new algorithms as well as easily reproducible results.

\subsection{Implementation perspective}

Algorithm development is an iterative process of designing and analyzing an
algorithm, prototyping and benchmarking it on small-sized test problems, and
making modifications until it performs well. Successful designs are then used to
solve larger problems. This workflow typically involves implementing and testing
the algorithm on different computing environments such as, \eg, multi-core
personal computers, dedicated servers, on-demand cloud-computing nodes, and even
a cluster of low-power IoT-enabled devices.

From an implementer's perspective, we would like to have a flexible library that
allows us to implement growing number of algorithm primitives on a diversity of
computing architectures. For this reason, we aim for building a low-level
library with small code-print, high performance and as few dependencies as
possible that would also support high-level language integrations and
cloud-based solutions.

\section{Design}\label{sec:design}

Our observations on different platforms and algorithms in the previous sections
reveal a trade-off in algorithm design. Decomposing algorithms into their small,
essential policies facilitates behavior isolation and code reuse at the expense
of combinatorially increasing design choices (\cf\ Figure~\ref{fig:proxgrad}).
To facilitate code reuse while handling the combinatorial explosion of choices,
\texttt{POLO} relies on \emph{policy-based design} ~\cite{2005-Alexandrescu},
which combines \emph{multiple inheritance} and \emph{template programming} in a
clever way.

To visualize how policy-based design can be used to implement optimization
algorithms, let us first assume that we would like to compute functions $y =
\func[\bar{f}]{\bar{x}}$ of the average of a streaming sequence of vectors
$\{x_{k}\}$. In languages that support both multiple inheritance and template
programming such as C++, one can define this functionality in terms of its
primitives $\func[\bar{f}]{\cdot}$ and $\func[ave]{\cdot}$ as given in
Listing~\ref{lst:pollib}.
\begin{lstlisting}[caption = {Simple example that visualizes how policy-based
design can be used to support functions operating on the average of a streaming
vector sequence. Common methods for construction, initialization and
destruction are omitted for brevity.}, label = lst:pollib, escapechar = |]
/* Template class for f */
template <class value_t, template <class> class fbar,
          template <class> class ave>
struct f : fbar<value_t>, ave<value_t> {
  value_t operator()(const vector<value_t> &x) {
    vector<value_t> temp(x.size());
    ave<value_t>::operator()(begin(x), end(x), begin(temp));
    return fbar<value_t>::operator()(begin(temp), end(temp));
  }
};

/* Exponential Moving Average */
template <class value_t> struct ema {
  template <class InputIt, class OutputIt>
  OutputIt operator()(InputIt xb, InputIt xe, OutputIt tb) {
    size_t idx{0};
    while (xb != xe) {
      average[idx] = (1 - alpha) * average[idx] + alpha * *xb++;
      *tb++ = average[idx++];
    }
    return tb;
  }

  void set_parameters(const value_t alpha) { this->alpha = alpha; }
private:
  value_t alpha;
  vector<value_t> average;
};

/* Sum of squared elements of a vector */
template <class value_t> struct sumsquared {
  template <class InputIt>
  value_t operator()(InputIt xb, InputIt xe) {
    value_t result{0};
    while (xb != xe) {
      result += *xb * *xb;
      xb++;
    }
    return result;
  }
};
\end{lstlisting}
Above, the library code first defines a template class for the function
\texttt{f}, which publicly inherits from its generic template parameters
\texttt{fbar} and \texttt{ave}, both of which are templated against the value
type (\texttt{value\_t}) they are working on. The library implements one
concrete policy for each function primitive. As can be observed, \texttt{f} can
work on different precision value-types such as \texttt{float}, \texttt{double}
and \texttt{fixed}, and its policy classes can have internal states
(\texttt{average}) and parameters (\texttt{alpha}) that need be persistent and
modifiable. Note how \texttt{f} is merely a shell that integrates its policy
classes and determines \emph{how} they should be applied to achieve the desired
transformation while not knowing their internal details.

A user of the library can then mix and match the provided policy classes to
achieve the desired functionality. The compiler assembles only the functionality
which has been requested and optimizes the code for the specified value types.
If the user wants to have a different averaging function such as the cumulative
moving average or they want to get the maximum absolute value of the averaged
$x$, which are not originally provided by the library, they can simply implement
these policies and use them together with the rest (see
Listing~\ref{lst:poluser}).
\begin{lstlisting}[caption = {Simple example that visualizes how policy-based
design can help users pick and implement different policies in a library.},
label = lst:poluser]
/* Use library-provided policies */
f<float, sumsquared, ema> f1;
f1.set_parameters(0.5);
f1.initialize({0,0});
f1({3,4}); /* returns 6.25 */
f1({6,8}); /* returns 39.0625 */

/* Implement cumulative moving average */
template <class value_t> struct cma {
  template <class InputIt, class OutputIt>
  OutputIt operator()(InputIt xb, InputIt xe, OutputIt tb) {
    size_t idx{0};
    while (xb != xe) {
      average[idx] = (N * average[idx] + *xb++) / (N + 1);
      *tb++ = average[idx++];
    }
    N++;
    return tb;
  }

private:
  size_t N{0};
  vector<value_t> average;
};

/* Implement maximum absolute value */
template <class value_t> struct maxabs {
  template <class InputIt>
  value_t operator()(InputIt xb, InputIt xe) {
    value_t result{0};
    while (xb != xe) {
      value_t current = abs(*xb++);
      if (current > result)
        result = current;
    }
    return result;
  }
};

/* Use user-defined policies */
f<double, maxabs, cma> f2;
f2.initialize({0,0});
f2({3,4}); /* cumulative average becomes (3,4), f2 returns 4 */
f2({6,8}); /* cumulative average becomes (4.5,6), f2 returns 6 */
\end{lstlisting}

Based on these simple examples, we immediately realize the strengths of the
policy-based design technique. First, the library designer and the user only
need to write code for the distinct functionalities, and can then combine them
into exponentially many design combinations. Second, the library can support
user-defined functionalities, which are not known in advance, and the compiler
can optimize the executable code for the selected data types thanks to template
programming. Last, the library can also support some enriched functionality that
the library designer has not foreseen. For example, the library designer
provides a shell \texttt{f} that only enforces the call function
(\texttt{operator()}) on its policies. However, thanks to the inheritance
mechanism, the user can use some added functionality (\texttt{set\_parameters})
when the policy classes provide them. In short, policy-based design helps handle
combinatorially increasing design choices with linearly many policies while also
allowing for extensions and enriched functionalities in ways library designers
might not have predicted.

Optimization algorithms have a lot in common with the examples provided in
Listings~\ref{lst:pollib}~and~\ref{lst:poluser}. Any family of algorithms
defines, at each iteration, how a sequence of different transformations should
be applied to some shared variable until a stopping condition is met. However,
optimization algorithms are much more involved than the examples above. They
have state loggers, different samplers and termination conditions, as well as
different ways of parallelizing various parts of the code. Moreover, not all the
policies of optimization algorithms are compatible with each other. As a result,
a similar approach is needed to represent these algorithms in a flexible way
while also preventing incompatible policies. To this end, we follow the standard
C++ parallel algorithms library proposals~\cite{2013-Hoberock,2015-Hoberock} in
\texttt{POLO} to add \emph{execution} support to optimization algorithms, and
use \emph{type traits} to disable incompatible policies at compile-time.

In the rest of the section, we will revisit proximal gradient methods to show
our design choice, and list the executors provided in \texttt{POLO}.

\subsection{Revisiting Proximal Gradient Methods}

In \texttt{POLO}, we provide the template class \texttt{proxgradient} to
represent the family of proximal gradient methods. It consists of
\texttt{boosting}, \texttt{step}, \texttt{smoothing}, \texttt{prox} and
\texttt{execution} policies, all of which are templated against the real
value-type, such as \texttt{float} or \texttt{double}, and the integral
index-type, such as \texttt{int32\_t} or \texttt{uint64\_t} (see
Listing~\ref{lst:proxgrad} for an excerpt).
\begin{lstlisting}[caption = {\texttt{proxgradient} algorithm template
implemented in \texttt{POLO}. Its policy templates default to implement the
vanilla gradient-descent algorithm running serially without any proximal step.},
label = lst:proxgrad]
/* namespace polo::algorithm */
template <class value_t = double, class index_t = int,
      template <class, class> class boosting = boosting::none,
      template <class, class> class step = step::constant,
      template <class, class> class smoothing = smoothing::none,
      template <class, class> class prox = prox::none,
      template <class, class> class execution = execution::serial>
struct proxgradient : public boosting<value_t, index_t>,
                      public step<value_t, index_t>,
                      public smoothing<value_t, index_t>,
                      public prox<value_t, index_t>,
                      public execution<value_t, index_t> {
  /* constructors and initializers */

  template <class Loss, class Terminator, class Logger>
  void solve(Loss &&loss, Terminator &&terminator, Logger &&logger)
  {
    execution<value_t, index_t>::solve(this,
      std::forward<Loss>(loss),
      std::forward<Terminator>(terminator),
      std::forward<Logger>(logger)
    );
  }

  /* other member functions */
};
\end{lstlisting}

As can be observed, \texttt{proxgradient} is simply a shell that glues together
its essential policies by using policy-based design, and delegates how to solve
the optimization problem to its \texttt{execution} policy. Based on this excerpt
only, we can see that any policy class that implements the corresponding
\texttt{solve} member function can be an executor for the \texttt{proxgradient}
family. For example, one can write a simple \texttt{serial} executor as in
Listing~\ref{lst:serial}.
\begin{lstlisting}[caption = {\texttt{serial} execution policy class implemented
for the \texttt{proxgradient} family in \texttt{POLO}. Due to space
considerations, only the relevant lines are shown.}, label = lst:serial,
escapechar = |]
/* namespace polo::execution */
template <class value_t, class index_t> struct serial {
  /* defaulted constructors and assignments */

protected:
  /* initializers and other functions */

  template <class Algorithm, class Loss, class Terminator,
            class Logger>
  void solve(Algorithm *alg, Loss &&loss, Terminator &&terminate,
             Logger &&logger) {
    while (!std::forward<Terminator>(terminate)(k, fval, xb_c,
                                                xe_c, gb_c)) {
      fval = std::forward<Loss>(loss)(xb_c, gb);
      iterate(alg, std::forward<Logger>(logger));
    }
  }

private:
  template <class Algorithm, class Logger>
  void iterate(Algorithm *alg, Logger &&logger) {
    alg->boost(index_t(0), k, k, gb_c, ge_c, gb);|\label{lst:startline}|
    alg->smooth(k, k, xb_c, xe_c, gb_c, gb);
    const value_t step = alg->step(k, k, fval, xb_c, xe_c, gb_c);
    alg->prox(step, xb_c, xe_c, gb_c, xb);|\label{lst:endline}|
    std::forward<Logger>(logger)(k, fval, xb_c, xe_c, gb_c);
    k++;
  }

  index_t k{1};
  value_t fval{0};
  std::vector<value_t> x, g;
  value_t *xb; /* points to the first element of x */
  const value_t *xb_c; /* points to the first element of x */
  const value_t *xe_c; /* points to past the last element of x */
  value_t *gb; /* points to the first element of g */
  const value_t *gb_c; /* points to the first element of g */
  const value_t *ge_c; /* points to past the last element of g */
};
\end{lstlisting}

The \texttt{serial} execution policy implements, between lines
\ref{lst:startline}--\ref{lst:endline} in Listing~\ref{lst:serial}, the
pseudocode given between lines \ref{alg:startline}--\ref{alg:endline} in
Algorithm~\ref{alg:prox-grad}. In fact, it only keeps track of the shared
variables \texttt{k}, \texttt{fval}, \texttt{x} and \texttt{g}, and determines
the order of executions, while delegating the actual task to the policy classes.
In other words, it calls, sequentially in the given order, the \texttt{boost},
\texttt{smooth}, \texttt{step} and \texttt{prox} functions, which are inherited
in \texttt{proxgradient} from \texttt{boosting}, \texttt{smoothing},
\texttt{step} and \texttt{prox} policies, respectively, without knowing what
they actually do to manipulate $x$ and $g$. This orthogonal decomposition of
algorithms into policies makes it easier for both library writers and users to
implement these primitives independently of each other.

Next, let us see how a user of \texttt{POLO} can build their algorithms by
selecting from predefined policies in Listing~\ref{lst:usercode}.
\begin{lstlisting}[caption = {Example code that uses \texttt{POLO}'s policy
classes to assemble different proximal gradient methods.}, label = lst:usercode]
/* include libraries */
#include "polo/polo.hpp"
using namespace polo;
using namespace polo::algorithm;

/* Assemble Heavyball */
proxgradient<double, int, boosting::momentum,
             step::constant, smoothing::none> heavyball;

/* Assemble AdaGrad */
proxgradient<double, int, boosting::none,
             step::constant, smoothing::adagrad> adagrad;

/* Assemble Adam */
proxgradient<double, int, boosting::momentum,
             step::constant, smoothing::rmsprop> adam;

/* Assemble Nadam */
proxgradient<double, int, boosting::nesterov,
             step::constant, smoothing::rmsprop> nadam;
\end{lstlisting}
Here, the user can easily assemble their favorite proximal gradient methods,
such as Heavyball~\cite{1964-Polyak}, AdaGrad~\cite{2011-Duchi},
Adam~\cite{2014-Kingma} and Nadam~\cite{2016-Dozat}, by choosing among the
provided policy classes. Similarly, if the user would like to try a different
\texttt{smoothing} policy instead of \texttt{adagrad} and \texttt{rmsprop}, they
can define their \texttt{custom} smoother
(lines~\ref{lst:customdefb}--\ref{lst:customdefe} in
Listing~\ref{lst:customcode}), and later, use it with other policy classes
(lines~\ref{lst:customuseb}--\ref{lst:customusee}). In
Listing~\ref{lst:customcode}, the \texttt{custom} smoother is indeed the
\texttt{amsgrad} smoother, which improves the convergence properties of
Adam-like algorithms~\cite{2018-Reddi}.
\begin{lstlisting}[caption = {Example code that defines a custom policy class to
be used together with provided functionality.}, label = lst:customcode,
escapechar = |]
/* Define a custom smoother */
template <class value_t, class index_t> struct custom {|\label{lst:customdefb}|
  custom(const value_t beta = 0.99, const value_t epsilon = 1E-6)
      : beta{beta}, epsilon{epsilon} {}

  /* defaulted copy/move operations */

  template <class InputIt1, class InputIt2, class OutputIt>
  OutputIt smooth(const index_t, const index_t, InputIt1 xbegin,
                  InputIt1 xend, InputIt2 gprev, OutputIt gcurr) {
    value_t g_val{0};
    index_t idx{0};
    while (xbegin != xend) {
      xbegin++;
      g_val = *gprev++;
      nu[idx] = beta * nu[idx] + (1 - beta) * g_val * g_val;
      nu_hat[idx] = std::max(nu_hat[idx], nu[idx]);
      *gcurr++ = g_val / (std::sqrt(nu_hat[idx]) + epsilon);
      idx++;
    }
    return gcurr;
  }

protected:
  void parameters(const value_t beta, const value_t epsilon) {
    this->beta = beta;
    this->epsilon = epsilon;
  }

  template <class InputIt> void initialize(InputIt xbegin,
                                           InputIt xend) {
    nu = std::vector<value_t>(std::distance(xbegin, xend));
    nu_hat = std::vector<value_t>(nu);
  }

  ~custom() = default;

private:
  value_t beta{0.99}, epsilon{1E-6};
  std::vector<value_t> nu, nu_hat;
};|\label{lst:customdefe}|

/* Assemble an algorithm using the custom smoother */
proxgradient<double, int, boosting::momentum,|\label{lst:customuseb}|
             step::constant, custom> myalg;|\label{lst:customusee}|
\end{lstlisting}

Finally, it is also worth mentioning that, in Listing~\ref{lst:usercode}, the
unused \texttt{prox} and \texttt{execution} policies default to \texttt{none}
and \texttt{serial}, respectively, which means that the proximal operator is
identity and the algorithm runs sequentially on one CPU. By mixing in different
choices of these policies, the user can extend the respective algorithms'
capabilities to handle different proximal terms on different, supported
executors.

\subsection{Provided Executors}\label{sub:executors}

\begin{figure}
  \centering
  \includegraphics{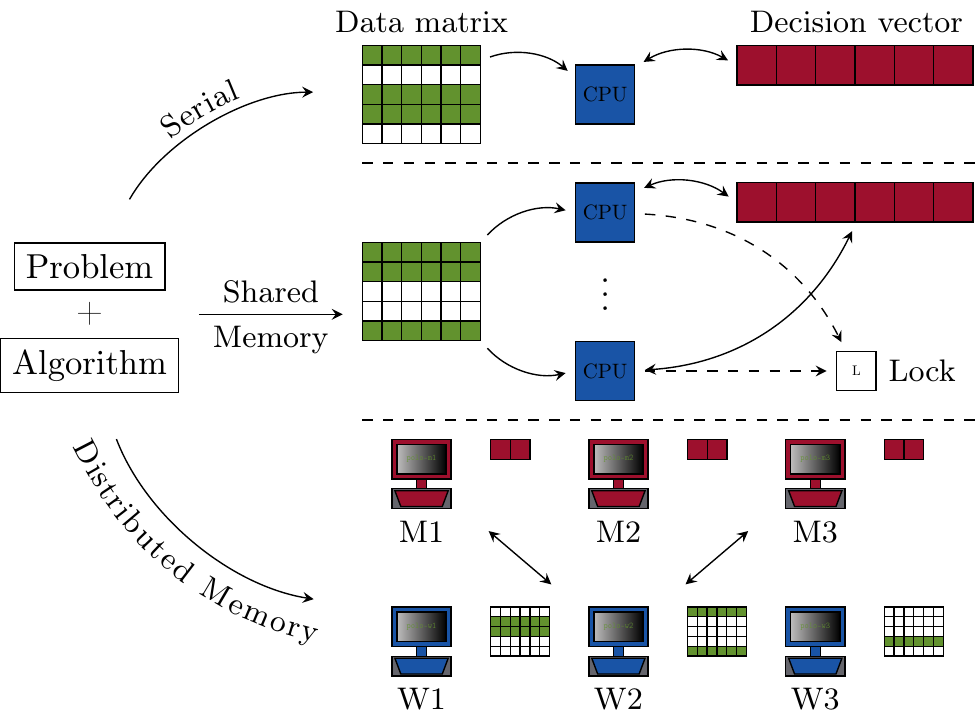}
  \caption{Supported platforms in \texttt{POLO} for the \texttt{proxgradient}
    family. When the problem data can fit on a single computer, the user can
    select among \texttt{serial}, \texttt{consistent} and \texttt{inconsistent}
    executions. The \texttt{consistent} and \texttt{inconsistent} executions
    differ from each other in that the former uses mutexes to lock (L) the whole
    decision vector when accessing the coordinates whereas the latter uses
    atomic operations on individual elements in shared-memory platforms. The
    library also supports distributed-memory executions when the problem data
    and decision vector are shared among different worker and master
    machines.}\label{fig:executors}
\end{figure}

\texttt{POLO} provides 4 different \texttt{execution} policy classes for the
\texttt{proxgradient} family to support 3 major computing platforms (see
Figure~\ref{fig:executors}):
\begin{enumerate}
  \item \texttt{serial} executor for sequential runs,
  \item \texttt{consistent} executor, which uses mutexes to lock the whole
    decision vector for consistent reads and writes, for shared-memory
    parallelism,
  \item \texttt{inconsistent} executor, which uses atomic operations to allow
    for inconsistent reads and writes to individual coordinates of the decision
    vector, for shared-memory parallelism, and,
  \item \texttt{paramserver} executor, which is a lightweight implementation of
    the Parameter Server~\cite{2013-Li} similar to that discussed
    in~\cite{2017-Xiao}, for distributed-memory parallelism.
\end{enumerate}

When the problem data can fit in a single memory space, users of \texttt{POLO}
can choose between the provided \texttt{serial}, \texttt{consistent} and
\texttt{inconsistent} execution policy classes. As we have partly seen in
Listing~\ref{lst:serial}, the \texttt{serial} execution simply implements the
pseudocode in Algorithm~\ref{alg:prox-grad}, and runs on a single CPU. When
multiple CPUs are available to speed-up the computations, users can choose
between \texttt{consistent} and \texttt{inconsistent} execution policies. In
both executors, updates to the decision vector are all well-defined. The name
\texttt{inconsistent} comes from the fact that it uses custom floating-point
types that support atomic operations, until these operations are supported
officially with C++20, inside the decision vector. Accessing and modifying
individual coordinates of the decision vector simultaneously from different
worker processes without locking them usually results in faster convergence at
the expense of suboptimal results, which can be controlled in certain
algorithm-problem pairs, especially when updates are relatively sparse.

\begin{figure}
  \centering
  \includegraphics{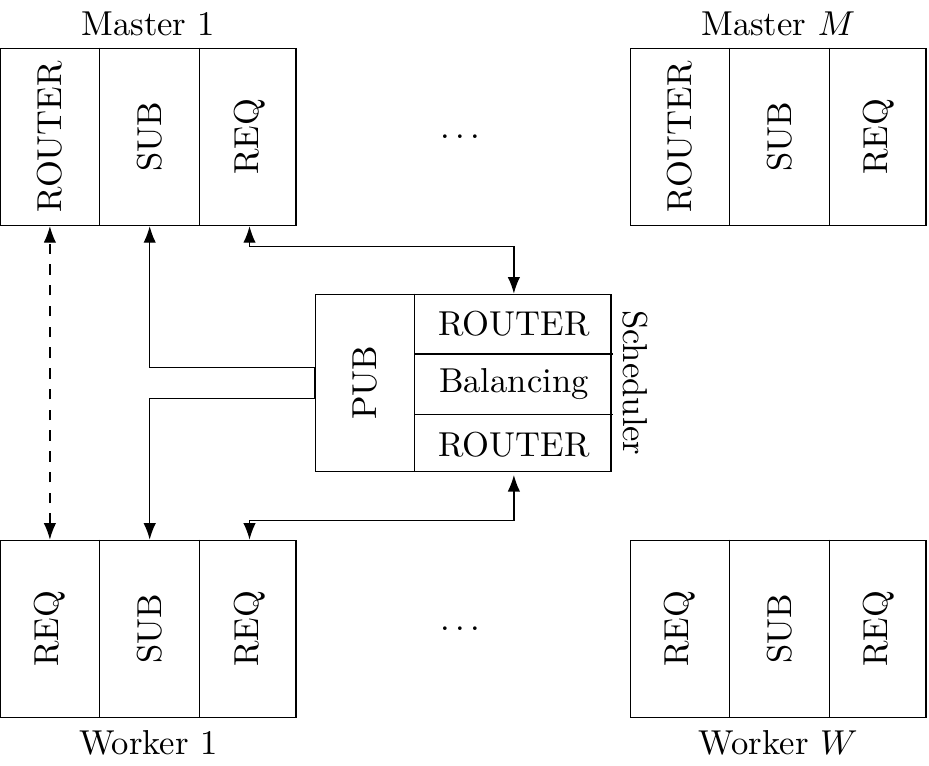}
  \caption{Implementation of Parameter Server in \texttt{POLO}. Arrow directions
    represent the flow of information. Solid lines represent static connections
    that have to be established at all times, whereas dashed lines represent
    dynamic connections that are established when
    needed.}\label{fig:paramserver}
\end{figure}

\texttt{POLO} also provides \texttt{paramserver} execution policy to support
computations which involve problem data scattered among different computing
nodes. Such solutions are needed when the data is either too large to fit on a
single node or it has to be kept at a certain place due to storing and accessing
requirements. The \texttt{paramserver} executor is not a full-fledged Parameter
Server~\cite{2013-Li} implementation but rather a lightweight distributed memory
executor similar to that mentioned in DSCOVR~\cite{2017-Xiao}\footnote{At the
time of writing this text, the Parameter Server website was not functioning and
there was no available code for either the Parameter Server or DSCOVR.}. It uses
\texttt{cURL}~\cite{2018-curl} and \texttt{ZMQ}~\cite{2017-ZMQ} libraries for
message passing, and \texttt{cereal}~\cite{2018-USCiLab} for serialization of
the transmitted data in a portable way. The \texttt{paramserver} executor has 3
main agents (see Figure~\ref{fig:paramserver}). The \texttt{scheduler} is the
only static node in the network. It is responsible for bookkeeping of connected
\texttt{master} nodes, publishing (PUB) global messages to the subscribed (SUB)
\texttt{master} and \texttt{worker} nodes, and directing \texttt{worker} nodes
to the corresponding \texttt{master} nodes when \texttt{worker} nodes need to
operate on specific parts of the decision vector. The \texttt{master} nodes
share the full decision vector in a linearly load-balanced fashion. They are
responsible for receiving gradient updates from \texttt{worker} nodes, taking
the \texttt{prox} step and sending the updated decision vector to the
\texttt{worker} nodes when they request (REQ). Finally, \texttt{worker} nodes
share the smooth loss-function data, and they can join the network dynamically.
Based on the sampling they use, \texttt{worker} nodes request a list of
\texttt{master} nodes that keep the corresponding part of the decision vector,
and they establish connections with them to communicate decision vectors and
(partial) gradient updates.

\section{Examples}\label{sec:examples}

In this section, we will demonstrate how \texttt{POLO} helps us implement
state-of-the-art optimization algorithms easily on different executors To this
end, we will first solve a randomly generated~\cite{1984-Lenard} unconstrained
quadratic problem using different serial algorithms from the proximal gradient
family. Then, we will solve a regularized logistic regression problem on the
\texttt{rcv1}~\cite{2004-Lewis} dataset using different algorithms in various
ways of parallelism.

We provide a sample code listing (Listing~\ref{lst:samplelisting}) in the
appendix for the examples we cover in this section, and state which parts of the
listing that need to be changed for each example.

\subsection{Unconstrained Quadratic Programming}

\begin{figure}[ht]
  \centering
  \begin{tabular}{@{}cc@{}}
    \includegraphics{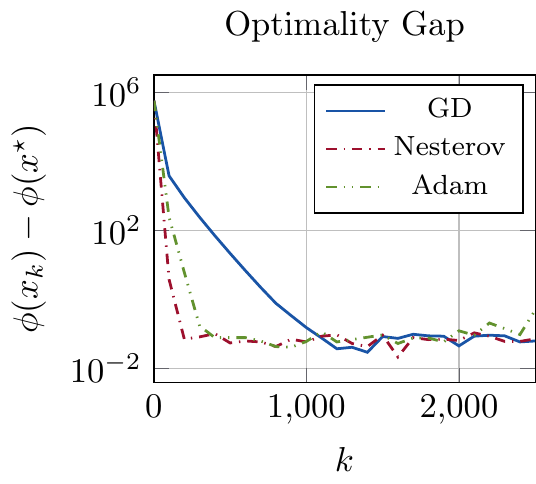} &
    \includegraphics{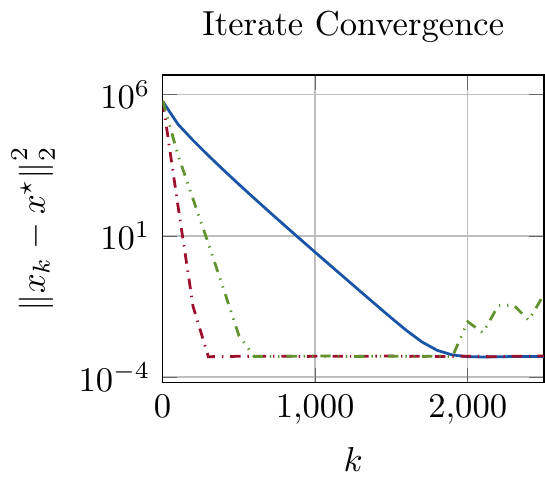}
  \end{tabular}
  \caption{Performance comparisons of different algorithms on the randomly
    generated QP problem. Adding \texttt{boosting} to full gradient iterations
    increases the rate of convergence, whereas adding \texttt{smoothing} results
    in a slightly slower convergence rate than that of \texttt{boosting}-only
    methods.}\label{fig:qp-example}
\end{figure}

In \texttt{POLO}, we provide a QP generator that follows the approach described
in~\cite{1984-Lenard} to create randomly generated QP problems on the form:
\begin{align*}
  \begin{aligned}
    & \minimize_{x \in \set{R}^{d}}
    & & \norm{\tilde{Q}x - \tilde{q}}_{2}^{2} \\
    & \text{subject to}
    & & Ax \leq b \,,
  \end{aligned}
\end{align*}
where $A$ is an $m \times d$ matrix, $b$ is an $m$-vector, $\tilde{Q}$ is an $n
\times d$ matrix, and $\tilde{q}$ is an $n$-vector. We use the generator to
create an unconstrained QP problem
\begin{align*}
  \begin{aligned}
    & \minimize_{x \in \set{R}^{d}}
    & & \frac{1}{2} x^{\top} Q x + q^{\top} x
  \end{aligned}
\end{align*}
with $d = 10000$, $\mu = \lambda_{\text{min}}(Q) = 1/20$ and $L =
\lambda_{\text{max}}(Q) = 20$. We use vanilla gradient descent, Nesterov
momentum and Adam to solve the problem. For this, one needs to change lines
\ref{lst:samplealgb}--\ref{lst:samplealge} of Listing~\ref{lst:samplelisting} to
have the following:
\begin{lstlisting}
#ifdef GD
algorithm::proxgradient<value_t, index_t> alg;
alg.step_parameters(2 / (mu + L));
alg.initialize(x0);
#elif defined NESTEROV
algorithm::proxgradient<value_t, index_t, boosting::nesterov> alg;
alg.boosting_parameters(0.9, 1 / L);
alg.initialize(x0);
#else
algorithm::proxgradient<value_t, index_t, boosting::momentum,
                        step::constant, smoothing::rmsprop> alg;
alg.step_parameters(0.08);
alg.boosting_parameters(0.9, 0.1);
alg.smoothing_parameters(0.999, 1E-8);
alg.initialize(x0);
#endif
\end{lstlisting}
In the end, we compile the code with appropriate definitions, and run the code
until the termination criterion is satisfied. We post-process the iterates
generated by the algorithms, and present the optimality gap and iterate
convergence in Figure~\ref{fig:qp-example}.

\subsection{Regularized Logistic Regression}

In \texttt{POLO}, we also provide various convenience functions such as dataset
readers, loss abstractions, samplers and basic matrix algebra functionality that
uses the installed BLAS/LAPACK implementations. It is worth noting that
\texttt{POLO} is not designed as a replacement for mature linear algebra
packages such as \texttt{Eigen}~\cite{2010-Guennebaud} and
\texttt{Armadillo}~\cite{2016-Sanderson}, and it can be used in combination with
these libraries when working with matrices.

In this subsection, we use the convenience functions to define a regularized
logistic loss
\begin{align*}
  \begin{aligned}
    & \minimize_{x \in \set{R}^{d}}
    & & \sum_{n = 1}^{N} \log(1 + \exp(-b_{i}\dotprod{a_{i}}{x})) +
      \lambda_{1} \norm{x}_{1}^{1} \,,
  \end{aligned}
\end{align*}
where the pair $\lbrace a_{i}, b_{i} \rbrace$ is the feature vector and the
corresponding label, respectively, of each sample in a given dataset, while
$\lambda_{1}$ is the regularization parameter. For this example, we will read
the \texttt{rcv1}~\cite{2004-Lewis} dataset in LIBSVM format. The dataset is
sparse with $N = 697641$ samples and $d = 47236$ features.

\begin{figure}[ht]
  \centering
  \begin{tabular}{@{}cc@{}}
    \includegraphics{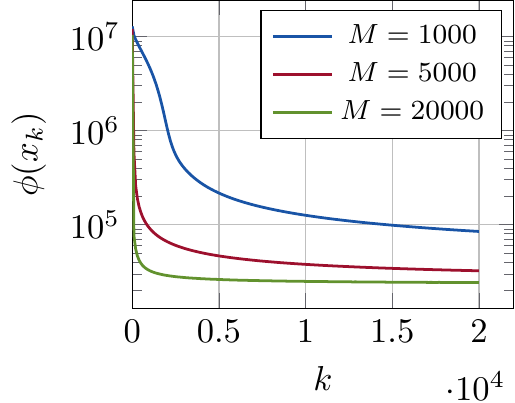} &
    \includegraphics{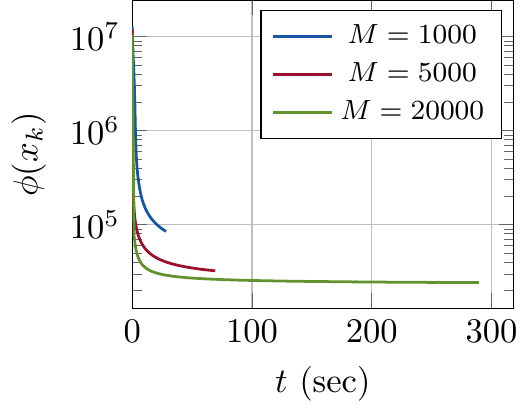}
  \end{tabular}
  \caption{Total function losses for AMSGrad with respect to the iteration count
    $k$ (left) and wall-clock time $t$ (right) on the \texttt{rcv1} dataset. As
    the number of mini-batches $M$ increases, the convergence rate (in terms of
    iteration count) increases and the error decreases, as expected. Moreover,
    the total computation time increases with the increasing number of
    mini-batches.}\label{fig:logloss-serial}
\end{figure}

We first assemble a different serial algorithm, AMSGrad~\cite{2018-Reddi}, and
use mini-batch gradients as the gradient surrogate, where mini-batches are
created by sampling the component functions uniformly at random. Necessary
changes to the sample code in Listing~\ref{lst:samplelisting} are as follows:
\begin{lstlisting}[escapechar = |]
/* Remove auxiliary variables between lines |\textcolor{green}{\ref{lst:sampleauxb}--\ref{lst:sampleauxe}}| in Listing |\textcolor{green}{\ref{lst:samplelisting}}| */
/* Replace the loss definition on line |\textcolor{green}{\ref{lst:sampleloss}}| in Listing |\textcolor{green}{\ref{lst:samplelisting}}| */
const string dsfile = "rcv1"; /* name of the dataset file */
auto dataset = utility::reader<value_t, index_t>::svm({dsfile});
loss::logistic<value_t, index_t> logloss(dataset);

const index_t K = 20000;
const index_t M = 1000; /* mini-batch size */
const index_t N = dataset.nsamples();
const index_t d = dataset.nfeatures();
const value_t L = 0.25 * M; /* L_i = 0.25; rcv1 is normalized */
const index_t B = N / M; /* number of mini-batches */
const value_t lambda1 = 1e-4;
/* Replace the algorithm between lines |\textcolor{green}{\ref{lst:samplealgb}--\ref{lst:samplealge}}| in Listing |\textcolor{green}{\ref{lst:samplelisting}}| */
algorithm::proxgradient<value_t, index_t, boosting::momentum,
                        step::constant, smoothing::amsgrad,
                        prox::l1norm> alg;

alg.step_parameters(1. / B); /* tuned for reasonable performance */
alg.boosting_parameters(0.9, 0.1); /* default suggested */
alg.smoothing_parameters(0.999, 1E-8); /* default suggested */
alg.prox_parameters(lambda1);
auto loss = [&](const value_t *x, value_t *g,
                const index_t *ibegin,
                const index_t *iend) -> value_t {
  return logloss.incremental(x, g, ibegin, iend);
};
/* Finally, change the solve method on line |\textcolor{green}{\ref{lst:samplesolve}}| in Listing |\textcolor{green}{\ref{lst:samplelisting}}| */
utility::sampler::uniform<index_t> sampler;
sampler.parameters(0, N - 1);
alg.solve(loss, utility::sampler::component, sampler, M,
          utility::terminator::maxiter<value_t, index_t>(K),
          logger);
\end{lstlisting}
As before, we compile and run the algorithm, this time with different mini-batch
sizes $M$, and then we reconstruct the total function loss values from the
traces saved in the \texttt{logger}. We report the results in
Figure~\ref{fig:logloss-serial}.

Next, we use the distributed memory functionality of \texttt{POLO} to solve
logistic regression on a parameter server architecture. To demonstrate how
\texttt{POLO} can be used equally well on embedded systems as on
high-performance computers, we first form a cluster of low-power devices. In the
cluster, we assign 5 Raspberry PI2 devices as the worker nodes, each of which
has parts of the samples in the \texttt{rcv1} dataset, and an NVIDIA Jetson TK1
as the single master, which only keeps track of the decision vector and does the
$\ell_{1}$-norm ball projection. We assign a laptop as the scheduler, which
orchestrates the communication between the static master nodes and the
dynamically joining worker nodes in the network. For this example, we choose to
use the proximal incremental aggregated gradient (PIAG)~\cite{2016-Aytekin}
method to solve the problem. Necessary modifications are given in the next
listing.
\begin{lstlisting}[escapechar = |]
/* Each worker node has their own local dataset */
/* Change the line |\textcolor{green}{\ref{lst:sampleloss}}| in Listing |\textcolor{green}{\ref{lst:samplelisting}}| */
#ifdef WORKER
const string dsfile = "rcv1";
auto dataset = utility::reader<value_t, index_t>::svm({dsfile});
loss::logistic<value_t, index_t> logloss(dataset);

const index_t N = dataset.nsamples();
const index_t d = dataset.nfeatures();
const value_t L = 0.25 * N;
#else
/* scheduler and master(s) need not know the loss */
auto loss = nullptr;
#endif
/* Change the lines |\textcolor{green}{\ref{lst:samplealgb}--\ref{lst:samplealge}}| in Listing |\textcolor{green}{\ref{lst:samplelisting}}| */
const value_t lambda1 = 1e-4;
const index_t K = 100000;

algorithm::proxgradient<value_t, index_t, boosting::aggregated,
                        step::constant, smoothing::none,
                        prox::l1norm,
                        execution::paramserver::executor> alg;

alg.step_parameters(1 / L);
alg.prox_parameters(lambda1);

execution::paramserver::options psopts;
/* workers timeout after 10 seconds of inactivity */
psopts.worker_timeout(10000);
/* master's own IP address; used by the worker */
psopts.master(maddress, 50000);
/* scheduler's address & ports; needed globally */
psopts.scheduler(saddress, 40000, 40001, 40002);

alg.execution_parameters(psopts);
\end{lstlisting}

\begin{figure}
  \centering
  \includegraphics{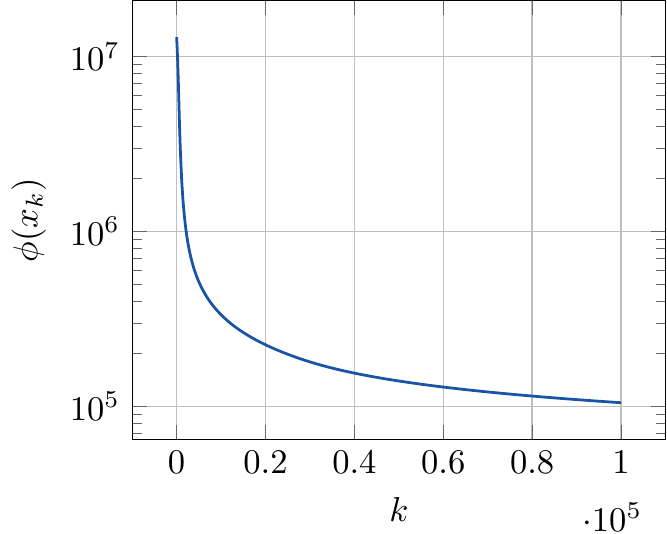}
  \caption{Total function loss on the \texttt{rcv1} dataset for PIAG executed on
  the Parameter Server executor provided in \texttt{POLO}.}
  \label{fig:logloss-ps}
\end{figure}

We compile the code three times, one for each of the worker, master and
scheduler agents in the distributed executor. We run the algorithm and
post-process the iterates logged by the master node. We report the result in
Figure~\ref{fig:logloss-ps}.

We finally recompiled the same code and ran it on virtual machines in a cloud
computing environment at Ericsson Research Data Center. No changes, apart from
reconfiguration of IP addresses, were needed to the code. Although the
wall-clock time performance improved significantly, the per-iterate convergence
is similar to the one for the PI-cluster and therefore not reported here.

\section{High-Level Language Integration}

Even though \texttt{POLO} helps us prototype ideas and implement
state-of-the-art optimization algorithms easily in C++, it is often even more
convenient to use a high-level language such as MATLAB, Python and Julia to
access data from different sources, define and automatically differentiate loss
functions, and visualize results. To be able to support integration with the
high-level languages, we also provide a C-API that implements the influential
algorithms listed in Table~\ref{tbl:proxgrad} on the executors covered in
Section~\ref{sub:executors}. Using the provided C header file and the compiled
library, users can access these algorithms from their favorite high-level
language.

Central to the C-API is a loss abstraction that allows for arbitrary loss
function implementations in any high-level language. The function signature of
the loss abstraction is given below.
\begin{lstlisting}[caption = {Loss abstraction used in the C-API. \texttt{data}
is an opaque pointer used to support callback functions from high-level
languages in the C library.}, escapechar = |, label = lst:lossabstraction]
template <class value_t>
using loss_t = value_t (*)(const value_t *x, value_t *g, void *data);
\end{lstlisting}
Any given loss should read $x$, write the gradient of $F$ at $x$ into $g$, and
return $\func[F]{x}$. The signature also includes an opaque pointer to
implementation specific data. For example, it could point to some loss function
object implemented in a high-level language, which calculates the loss and
gradient at $x$. Next, the C-API defines algorithm calls for a set of algorithm
types and executors as follows.
\begin{lstlisting}[escapechar = |]
using value_t = double; /* double-precision floating points */
void run_serial_alg(const value_t *xbegin, const value_t *xend,
                    loss_t<value_t> loss_fcn, void *loss_data) {
  serial_alg_t alg;
  alg.initialize(xbegin, xend);
  alg.solve([=](const value_t *xbegin, value_t *gbegin) {
    return loss_fcn(xbegin, gbegin, loss_data);
  });
}
\end{lstlisting}
The \texttt{serial\_alg\_t} is varied to specific algorithm types through
different policy combinations. For example, the following defines regular
gradient descent with a serial executor.
\begin{lstlisting}[escapechar = |]
using namespace polo;
using namespace polo::algorithm;
using serial_alg_t =
  proxgradient<double, int, boosting::none, step::constant,
               smoothing::none, prox::none, execution::serial>;
\end{lstlisting}
High-level languages can then make use of these algorithms by supplying loss
functions that adhere to the abstraction in Listing~\ref{lst:lossabstraction}.

The C-API also defines a set of custom policies within the \texttt{proxgradient}
algorithm family. These policies have empty implementations that can be filled
in by high-level languages. This feature is implemented using abstractions with
opaque pointers, similar to the loss abstraction in
Listing~\ref{lst:lossabstraction} and gives the user extended functionality
outside the range of precompiled algorithms in the C-API. To illustrate this
versatility, we have implemented \texttt{POLO.jl} in the Julia language.
\texttt{POLO.jl} includes Julia implementations of all \texttt{proxgradient}
policies available in \texttt{POLO}. Through the custom policy API, any
algorithm type can be instantiated and tested interactively. In addition, it
allows us to implement other high-level abstractions in Julia such as custom
termination criteria and advanced logging. As an example, the following
implements Adam in \texttt{POLO.jl}
\begin{lstlisting}[escapechar = |]
Adam(execution::ExecutionPolicy) =
  ProxGradient(execution, Boosting::Momentum(), Step.Constant(),
               Smoothing.RMSprop(), Prox.None())
\end{lstlisting}
In essence, \texttt{POLO.jl} provides a dynamic setting for algorithm design
while also enabling the use of the powerful executors implemented in
\texttt{POLO}. Due to current limitations in the Julia language, the
multi-threaded executors could not yet be applied.

\section{Conclusion and Future Work}\label{sec:conclusion}

In this paper, we have presented \texttt{POLO}, an open-source, header-only, C++
template library. Compared to other low-level generic optimization libraries,
\texttt{POLO} is focused on algorithm development. With its proper algorithm
abstractions, policy-based design approach and convenient utilities,
\texttt{POLO} not only offers optimization and machine-learning researchers an
environment to prototype their ideas flexibly without losing much from
performance but also allows them to test their ideas and other state-of-the-art
optimization algorithms on different executors easily. \texttt{POLO} is still a
work in progress, and we are working on adding other families of algorithms and
different executors to the library as well as wrapping the C-API as convenient
packages in high-level languages such as Python and Julia.

\section*{Acknowledgment}

This research is sponsored in part by the Swedish Research Council under project
grant ``Scalable optimization: dynamics, architectures and data dependence,''
and by the Swedish Foundation for Strategic Research under project grant
``Societal-scale cyber-physical transport systems.'' We also thank Ericsson
Research for their generosity in letting us use the computing resources at
Ericsson Research Data Center for our experiments.

\appendix
\section{Sample Code Listing}

\begin{lstlisting}[caption = {Sample code listing to reproduce experiments in
Section~\ref{sec:examples}. Only \texttt{POLO} related functionality is shown
for brevity.}, label = lst:samplelisting, escapechar = |]
/* Include necessary libraries */
using index_t = int32_t;
using value_t = float;

int main(int argc, char *argv[]) {
  /* Define auxiliary variables */
  const index_t d = 10000;|\label{lst:sampleauxb}|
  const index_t K = 2500;
  const value_t L = 20;
  const value_t mu = 1 / L;|\label{lst:sampleauxe}|

  /* Define a loss function */
  problem::qp<value_t, index_t> qp(d, L);|\label{lst:sampleloss}|

  /* Randomly initialize the starting point */
  mt19937 gen(random_device{}());
  normal_distribution<value_t> normal(5, 3);
  vector<value_t> x0(d);
  for (auto &val : x0)
    val = normal(gen);

  /* Select, configure and initialize an algorithm */
  algorithm::proxgradient<value_t, index_t> alg;|\label{lst:samplealgb}|
  alg.step_parameters(2 / (mu + L));
  alg.initialize(x0);|\label{lst:samplealge}|

  /* Use a decision logger */
  utility::logger::decision<value_t, index_t> logger;

  /* Terminate after K iterations */
  utility::terminator::maxiter<value_t, index_t> maxiter(K);

  /* Solve the problem */
  alg.solve(qp, maxiter, logger);|\label{lst:samplesolve}|

  /* Post-process the logged states, i.e., function values,
     decision vector, etc.
  */

  return 0;
}
\end{lstlisting}

\printbibliography

\end{document}